
\def\LaTeX{%
  \let\Begin\begin
  \let\End\end
  \let\salta\relax
  \let\finqui\relax
  \let\futuro\relax}

\def\UK{\def\our{our}\let\sz s}
\def\USA{\def\our{or}\let\sz z}

\UK 



\LaTeX

\USA


\salta

\documentclass[twoside,12pt]{article}
\setlength{\textheight}{24cm}
\setlength{\textwidth}{16cm}
\setlength{\oddsidemargin}{2mm}
\setlength{\evensidemargin}{2mm}
\setlength{\topmargin}{-15mm}
\parskip2mm


\usepackage[usenames,dvipsnames]{color}
\usepackage{amsmath}
\usepackage{amsthm}
\usepackage{amssymb,bbm}
\usepackage[mathcal]{euscript}

\usepackage{cite}

\usepackage{hyperref}
\usepackage[shortlabels]{enumitem}



%
%

%
 
\definecolor{viola}{rgb}{0.3,0,0.7}
\definecolor{ciclamino}{rgb}{0.5,0,0.5}
\definecolor{blu}{rgb}{0,0,0.7}
\definecolor{rosso}{rgb}{0.85,0,0}

\def\gianni #1{{\color{red}#1}}
\def\pier #1{{\color{blue}#1}} 

\def\gianni #1{#1}
\def\pier #1{#1}






%

\finqui

\def\Beq{\Begin{equation}}
\def\Eeq{\End{equation}}

\def\Bthm{\Begin{theorem}}
\def\Ethm{\End{theorem}}
\def\Blem{\Begin{lemma}}
\def\Elem{\End{lemma}}

\def\Brem{\Begin{remark}\rm}
\def\Erem{\End{remark}}

\def\Bdim{\Begin{proof}}
\def\Edim{\End{proof}}
\def\Bcenter{\Begin{center}}
\def\Ecenter{\End{center}}
\let\non\nonumber




\def\step #1 \par{\medskip\noindent{\bf #1.}\quad}
\def\jstep #1: \par {\vspace{2mm}\noindent\underline{\sc #1 :}\par\nobreak\vspace{1mm}\noindent}

\def\aand{\quad\hbox{and}\quad}
\def\Lip{Lip\-schitz}
\def\Holder{H\"older}

\def\lhs{left-hand side}
\def\rhs{right-hand side}




\def\multibold #1{\def\arg{#1}%
  \ifx\arg\pto \let\next\relax
  \else
  \def\next{\expandafter
    \def\csname #1#1\endcsname{{\boldsymbol #1}}%
    \multibold}%
  \fi \next}

\def\pto{.}

\def\multical #1{\def\arg{#1}%
  \ifx\arg\pto \let\next\relax
  \else
  \def\next{\expandafter
    \def\csname cal#1\endcsname{{\cal #1}}%
    \multical}%
  \fi \next}

\def\multigrass #1{\def\arg{#1}%
  \ifx\arg\pto \let\next\relax
  \else
  \def\next{\expandafter
    \def\csname grass#1\endcsname{{\mathbb #1}}%
    \multigrass}%
  \fi \next}


\def\multimathop #1 {\def\arg{#1}%
  \ifx\arg\pto \let\next\relax
  \else
  \def\next{\expandafter
    \def\csname #1\endcsname{\mathop{\rm #1}\nolimits}%
    \multimathop}%
  \fi \next}

\multibold
qweryuiopasdfghjklzxcvbnmQWERTYUIOPASDFGHJKLZXCVBNM.  

\multical
QWERTYUIOPASDFGHJKLZXCVBNM.

\multigrass
QWERTYUIOPASDFGHJKLZXCVBNM.

\multimathop
diag dist div dom mean meas sign supp .

\def\Span{\mathop{\rm span}\nolimits}


\def\accorpa #1#2{\eqref{#1}--\eqref{#2}}
\def\Accorpa #1#2 #3 {\gdef #1{\eqref{#2}--\eqref{#3}}%
  \wlog{}\wlog{\string #1 -> #2 - #3}\wlog{}}


\def\separa{\noalign{\allowbreak}}

\def\somma #1#2#3{\sum_{#1=#2}^{#3}}

\def\graffe #1{\mathopen\{#1\mathclose\}}
\def\<#1>{\mathopen\langle #1\mathclose\rangle}
\def\norma #1{\mathopen \| #1\mathclose \|}
\def\aeO{\checkmmode{a.e.\ in~$\Omega$}}
\def\aeQ{\checkmmode{a.e.\ in~$Q$}}

\def\aet{\checkmmode{a.e.\ in~$(0,T)$}}
\def\aat{\checkmmode{for a.e.\ $t\in(0,T)$}}

\let\hat\widehat
\def\cpto{\,\cdot\,}

\def\iot {\int_0^t}
\def\ioT {\int_0^T}
\def\intQt{\int_{Q_t}}
\def\intQ{\int_Q}
\def\iO{\int_\Omega}

\def\dt{\partial_t}
\def\dn{\partial_{\nn}}
\def\ddt{\frac d{dt}}

\def\Pbln{P$_\nn$}
\def\pbln{P$_n$}

\let\emb\hookrightarrow
\def\cpto{\,\cdot\,}

\def\checkmmode #1{\relax\ifmmode\hbox{#1}\else{#1}\fi}


\let\erre\grassR

\def\erren{\erre^n}




\def\genspazio #1#2#3#4#5{#1^{#2}(#5,#4;#3)}
\def\spazio #1#2#3{\genspazio {#1}{#2}{#3}T0}

\def\L {\spazio L}
\def\H {\spazio H}
\def\W {\spazio W}

\def\C #1#2{C^{#1}([0,T];#2)}


\def\Lx #1{L^{#1}(\Omega)}
\def\Hx #1{H^{#1}(\Omega)}
\def\Wx #1{W^{#1}(\Omega)}

\def\LQ #1{L^{#1}(Q)}

\def\CQ #1{C^{#1}(\overline Q)}

\def\Ldue{\Lx 2}
\def\Linfty{\Lx\infty}

\def\Huno{\Hx 1}
\def\Hdue{\Hx 2}



\let\eps\varepsilon

\let\phi\varphi

\let\TeXchi\chi                         
\newbox\chibox
\setbox0 \hbox{\mathsurround0pt $\TeXchi$}
\setbox\chibox \hbox{\raise\dp0 \box 0 }
\def\chi{\copy\chibox}



\def\Beta{\widehat\beta}
\def\ueps{u_\eps}

\def\ej{e_j}
\def\ei{e_i}
\def\lambdai{\lambda_i}
\def\lambdaj{\lambda_j}
\def\phii{\pier{\phi_{n,i}}}
\def\phij{\pier{\phi_{n,j}}}
\def\mui{\pier{\mu_{n,i}}}
\def\muj{\pier{\mu_{n,j}}}

\def\bphi{\boldsymbol\phi}  
\def\bmu{\boldsymbol\mu}

\def\Vn{V_n}
\def\phin{\phi_n}
\def\mun{\mu_n}

\def\vn{v_n}
\def\Pn{\grassP_n}  

\def\phit{\phi_\tau}
\def\mut{\mu_\tau}
\def\wt{w_\tau}

\def\phibar{\phi_\Omega}
\def\mubar{\mu_\Omega}
\def\phizbar{\phiz_\Omega}
\def\rhozbar{\rhoz_\Omega}
\def\gbar{g_\Omega}

\def\vbar{v_\Omega}
\def\zbar{z_\Omega}
\def\zetabar{\zeta_\Omega}

\def\phiz{\phi^0}
\def\rhoz{\rho^0}

\def\Vp{{V^*}}
\def\Wp{W^*}
\def\Zp{Z^*}

\def\soluz{(\phi,\mu,w)}
\def\soluzbis{(\phi,\mu)}
\def\soluzn{(\phin,\mun)}
\def\soluzt{(\phit,\mut,\wt)}

\def\normaV #1{\norma{#1}_V}
\def\normaW #1{\norma{#1}_W}
\def\normaZ #1{\norma{#1}_Z}
\def\normaVp #1{\norma{#1}_*}

\def\CO{C_\Omega}
\def\Cbeta{C_\beta}
\def\cdelta{c_\delta}


\usepackage{amsmath}
\DeclareFontFamily{U}{mathc}{}
\DeclareFontShape{U}{mathc}{m}{it}%
{<->s*[1.03] mathc10}{}

\DeclareMathAlphabet{\mathscr}{U}{mathc}{m}{it}

\Begin{document}


%
\title{Hyperbolic relaxation\\ of a sixth-order Cahn--Hilliard equation}
\author{}
\date{}
\maketitle
\Bcenter
\vskip-1.8cm
{\large\sc Pierluigi Colli$^{(1,2)}$}\\
{\normalsize e-mail: {\tt pierluigi.colli@unipv.it}}\\
[0.25cm]
{\large\sc Gianni Gilardi$^{(1,3)}$}\\
{\normalsize e-mail: {\tt gianni.gilardi@unipv.it}}\\
[.5cm]
$^{(1)}$
{\small Dipartimento di Matematica ``F. Casorati'', Universit\`a di Pavia}\\
{\small via Ferrata 5, I-27100 Pavia, Italy}\\
[0.25cm]
$^{(2)}$
{\small Research Associate at the IMATI -- C.N.R. Pavia}\\
[0.25cm]
$^{(3)}$
{\small Istituto Lombardo, Accademia di Scienze e Lettere, Milano}\\
\Ecenter
%
\Begin{abstract}
\noindent
This work explores the solvability of a sixth-order Cahn--Hilliard equation with an inertial term, which serves as a relaxation of a higher-order variant of the classical Cahn--Hilliard equation. The equation includes a source term that 
disrupts the conservation of the mean value of the order parameter. The incorporation of additional spatial derivatives allows the model to account for curvature effects, leading to a more precise representation of isothermal phase 
separation dynamics. We establish the existence of a weak solution for the associated initial and boundary value problem under the assumption that the double-well-type nonlinearity is globally defined. Additionally, we derive 
uniform stability estimates, which enable us to demonstrate that any family of solutions satisfying these estimates converges in a suitable topology to the 
unique solution of the limiting problem as the relaxation parameter approaches zero. Furthermore, we provide an error estimate for specific norms of the difference between solutions in terms of the relaxation parameter.
\\[2mm]
\\[2mm]
\End{abstract}

%

\pagestyle{myheadings}
\newcommand\testopari{\sc Colli -- Gilardi}
\newcommand\testodispari{\sc Hyperbolic relaxation of a sixth-order Cahn--Hilliard equation}
\markboth{\testopari}{\testodispari}
%

\section{Introduction}
\label{INTRO}
\setcounter{equation}{0}
In a recent study~\cite{CGSS6} conducted by the authors in collaboration with A. Signori and J. Sprekels, a sixth-order partial differential equation of the Cahn--Hilliard type was examined and analyzed. Specifically, by introducing the chemical potential $\mu$ and the auxiliary variable $w$  associated with the phase (or order) parameter $\phi$, the study~\cite{CGSS6} focuses on the resulting system of three equations, namely:
\begin{align}
  & \dt\phi - \Delta\mu + \sigma\phi 
  = g 
  \qquad && \text{in $Q$}
  \label{Iprimazero}
  \\
  & -\Delta w + f'(\phi) w + \nu w
  = \mu
  \quad && \text{in $Q$}
  \label{Iseconda}
  \\
  & - \Delta\phi + f(\phi)
  = w
  \quad && \text{in $Q$.}
  \label{Iterza}
\end{align}
No-flux boundary conditions and an initial condition have been added to these equations to complete the system:
\begin{align}
  & \dn\mu
  =  \dn w  
  = \dn\phi  
  = 0
  \quad && \text{on $\Sigma$}
  \label{Ibc}
  \\
  & \phi(0) = \phiz
  \quad && \text{in $\Omega$}.
  \label{Iic}
\end{align}
\Accorpa\Ipblzero Iprimazero Iic
Here and throughout the sequel, $Q$ and $\Sigma$ are defined~as
\Beq
  Q := \Omega\times(0,T)
  \aand
  \Sigma := \Gamma\times(0,T)
  \label{defQS}
\Eeq
where $\Omega$ is the domain in $\erre^3$ where the evolution takes place,
$\Gamma:=\partial\Omega$ is its boundary and $T$~stands for a given final time.
In the above equations, $f$~is the derivative of a regular double well potential~$F$,
while $\sigma$ and $\nu$ are given constants.
Moreover, $g$~in \eqref{Iprimazero} is a prescribed forcing term, $\dn$~in \eqref{Ibc} denotes the outward normal derivative
and $\phiz$ in \eqref{Iic} is a given initial datum.
A~prototype of $F$ is the {\it classical regular potential}, given~by
\Beq
  F(s) = \frac 14 (s^2-1)^2 ,
  \, \hbox{ and consequently } \, f(s) = s^3 - s , \quad \hbox{for $s\in\erre$}.
  \label{regpot}
\Eeq
The chemical potential $\mu$ in equation \eqref{Iseconda} corresponds to the first variation of the free energy $\cal E$. Similarly, $w$ is the first 
variational  derivative of the Ginzburg--Landau free energy ${\cal G}$. Specifically, it holds that $\mu= {\delta {\cal E}}/ {\delta \phi}$ and $w= {\delta {\cal G}}/{\delta \phi}$, where
\Beq
  \calE(v) =
 {\cal F}(\phi) + {\nu} \,{\cal G}(\phi)=\frac 12 \iO \bigl( -\Delta v + f(v) \bigr)^2
  + \nu \iO \Bigl( \frac 12 \, |\nabla v|^2 + F(v) \Bigr). 
  \label{Ienergy}
\Eeq
The system \eqref{Iprimazero}--\eqref{Iterza} can be regarded as a variant of the classical fourth-order Cahn--Hilliard equation~\cite{CH} (see the recent monograph~\cite{CHbook} for an overview) in terms of the unknowns $\phi$,
 $\mu$, and $w$. The original Cahn--Hilliard equation is widely used to model {isothermal} phase separation in binary
 mixtures. The variable $\phi$ represents the local proportion of one component within the binary material and serves as an order parameter. For simplicity, it is typically normalized so that the pure states correspond to  $\phi = \pm 1$, while the region 
${\{ -1 < \phi < 1\}} $ represents the diffuse interface, which forms in a tubular neighborhood of the interface with a small thickness parameter.
Unlike the classical model, the sixth-order Cahn--Hilliard equation (see \cite{PZ2, SW, M1}) extends it by incorporating curvature effects and higher-order variations of the order parameter. This extension enables a more accurate representation of complex material systems with intricate interface structures. Furthermore, it has potential applications in the design and optimization of materials with tailored microstructures and properties.

In the paper~\cite{CGSS6} it is proved that  the initial and boundary value problem~\Ipblzero\ has a unique solution in a proper functional framework; moreover, an optimal control problem
is investigated in terms of a distributed control, represented by $g$ in equation~\eqref{Iprimazero}, and a suitable cost functional. 

In this paper, we shift our focus to exploring the solvability of a hyperbolic relaxation of problem~\Ipblzero. Additionally, we delve into the convergence and asymptotic behavior of solutions, examining how they align with the original problem as the relaxation parameter approaches zero. Specifically, we modify equation~\eqref{Iprimazero}~by
\Beq
  \tau \dt^2\phi + \dt\phi - \Delta\mu + \sigma\phi 
  = g 
  \qquad \text{in $Q$}
  \label{Iprimatau}
\Eeq
where $\tau$ is a positive coefficient. Coupling \eqref{Iprimatau} with \eqref{Iseconda} and \eqref{Iterza} leads to a sixth-order equation with the inertial term $\tau \dt^2\phi$.

Some remarks on systems of this type are in order. Over the past two decades, there has been a notable increase in mathematical research on hyperbolic variants of Cahn--Hilliard systems that include inertial terms. The introduction of hyperbolic relaxation terms was first proposed by Galenko~\cite{Gal} and later revisited by Galenko and other physicists \cite{GalLeb, GalJou, CCJ} to model strongly non-equilibrium decompositions. These decompositions arise due to rapid solidification under supercooling into the spinodal region, a phenomenon observed in certain materials such as glasses. 

Regarding the mathematical analysis of hyperbolic versions of the Cahn--Hilliard equations, we refer to \cite{Bon1, GGMP, GGPM, GraPie, GPS, GSZ1, GSZ2, GSSZ, Seg}, as well as the more recent contributions \cite{Bon2, Chen, CMY, Dor, DMP, SZ, SK, XL}. These studies address various aspects, including the existence and regularity of solutions in 2D and 3D, the incorporation of additional mass terms in the first equation~\eqref{Iprimazero} (such as the Oono term $\sigma \phi - g$), asymptotic behavior and attractors, approximations and discretized systems, and equations defined in unbounded domains.

A noteworthy contribution is the study by Chen and Liu~\cite{ChenLiu}, which investigates a sixth-order Cahn--Hilliard equation with an inertial term to model phase transition dynamics in ternary oil-water-surfactant systems. This analysis is conducted in the one-dimensional case with periodic boundary conditions. Additionally, the recent work~\cite{CS12} examines a hyperbolic relaxation of the viscous Cahn--Hilliard system, where the relaxation term involves the second time derivative of the chemical potential rather than the phase variable.

Returning to \eqref{Iprimatau} coupled with \eqref{Iseconda}--\eqref{Iic},
we note that an additional initial condition must be imposed, namely
\Beq
  \dt\phi(0) = \rhoz 
  \label{Icauchydt}
\Eeq
where $\rhoz$ is given, to properly define a suitable initial and boundary value problem. We establish that this modified problem admits at least one solution (in a weak sense) that also satisfies precise stability estimates. Furthermore, we demonstrate that any family $\graffe\soluzt$ of solutions to the relaxed problem, satisfying the stability estimates, converges in an appropriate topology to the unique solution $\soluz$ to problem \Ipblzero\ as $\tau$ tends to zero. Additionally, we derive an error estimate in weaker norms for the difference between solutions, expressed in terms of the relaxation parameter, and establish a convergence rate of order~\gianni{$\tau^{1/2}$}. 

The paper is structured as follows. In the next section, we introduce our assumptions and notations, present our main results, and compile some useful tools employed throughout the work. The proofs of our existence and asymptotic results are provided in Sections~\ref{EXISTENCE} and~\ref{ASYMPTOTICS}, respectively.


\section{Notation, assumptions and results}
\label{STATEMENT}
\setcounter{equation}{0}

As for the set $\Omega$ mentioned in the Introduction,
we assume that it is a bounded and connected open set in $\erre^3$ with smooth boundary $\Gamma$
and we denote by $|\Omega|$ its Lebesgue measure.
For any Banach space $X$ the notations $\norma\cpto_X$, $X^*$, and $\< \cpto,\cpto >_X$ indicate 
the corresponding norm, its dual space, and the related duality pairing between $X^*$ and~$X$,
with some exceptions for the notation of the norm that are listed below.
Next, we~set
\begin{align}
  & H := \Ldue , \quad  
  V := \Huno, \quad
  W := \graffe{v\in\Hdue: \ \dn v=0 \hbox{ on $\,\Gamma$}} 
  \non
  \\
  & \aand
  Z := \graffe{v\in W: \ \Delta v\in W}.
 \label{defspazi}
\end{align}
The norm in the special case $H$ is indicated by~$\norma\cpto$ (instead of $\norma\cpto_H$) to~simplify the notation.
We also denote by $(\cpto,\cpto)$ the standard inner product of~$H$.
The symbol $\norma\cpto_\infty$ might denote the norm in each of the spaces 
$\Linfty$, $\LQ\infty$ (see \eqref{defQS}) and $L^\infty(0,T)$, if no confusion can arise.

\Brem
\label{Density}
Clearly, $V$ is dense in~$H$ and $W$~is dense in~$V$.
It turns out that $Z$ is dense in~$W$, as we briefly sketch.
Given $u\in W$, we approximate it by the unique solution to the problem
\Beq
  - \eps \Delta\ueps + \ueps = u
  \quad \hbox{in $\Omega$}
  \aand
  \dn\ueps = 0
  \quad \hbox{on $\Gamma$}
  \label{ueps}
\Eeq
where $\eps$ is a positive parameter.
Clearly $\ueps\in W$, whence also $\Delta\ueps\in W$.
This implies that $\ueps\in Z$ and that
\Beq
  \eps \Delta^2 \ueps - \Delta\ueps = - \Delta u \,.
  \label{uepsbis}
\Eeq
Then, by testing \eqref{ueps} by $\ueps$ and $-\Delta\ueps$ and \eqref{uepsbis} by $-\Delta\ueps$,
we easily obtain that
\Beq
  \norma\ueps \leq \norma u \,, \quad
  \norma{\nabla\ueps} \leq \norma{\nabla u}
  \aand
  \norma{\Delta\ueps} \leq \norma{\Delta u}.
  \non
\Eeq
Hence, $\ueps$ has a subsequence that converges (weakly in principle) in~$W$ as $\eps$ tends to zero
and one immediately sees that its limit is~$u$.
This shows the expected density.
By the way, the whole family $\graffe{\ueps}$ converges to $u$ and the convergence is even strong, by uniform convexity.
\Erem

By the above remark,
we can adopt the usual framework of Hilbert triplets obtained by standard identifications.
Namely, we have that
\begin{align}
  & \< z,v >_V = (z,v) , \quad
  \< z,v >_W = \< z,v >_V,
  \aand
  \< z,v >_Z = \< z,v >_W,
  \non
  \\
  & \quad \hbox{for every $z\in H,\Vp,\Wp$ and $v\in V,W,Z$, respectively}.
  \non
\end{align}
This allows us to use the simpler symbol $\<\cpto,\cpto>$ without any subscript for the above duality pairings
whenever no confusion can arise.
In conclusion, we have that
\Beq
  Z \emb W \emb V \emb H \emb \Vp \emb \Wp \emb \Zp
  \label{embeddings}
\Eeq
with dense and compact embeddings. 

Now, we list our assumptions on the structure of the system to be analyzed.
As in \cite{CGSS6}, we assume:
\begin{align}
  & \sigma, \, \lambda \in (0,+\infty) 
  \aand
  \nu \in \erre.
  \label{hpconstants}
  \\
  & F \in C^4(\erre) \ \hbox{can be written as} \quad
  F(s) = \Beta(s) - \frac \lambda 2 \, s^2 \,,
  \quad s \in \erre \,,
  \quad \hbox{with $\Beta$ convex} .
  \label{hpBeta}
\end{align}
We set
\Beq
  f := F' \aand \beta := {\Beta\,}'
  \label{defbeta}
\Eeq
so that 
\begin{align}
	f(s)= F'(s)= \beta (s) - \lambda s
	\aand 
	f'(s)= \beta'(s) - \lambda,
  \quad \hbox{for every $s\in\erre$\pier.}
  \label{calcolof}
\end{align}
\pier{Concerning $\beta$ we follow exactly the same framework as in \cite[assumptions~(2.6)--(2.8)]{CGSS6} and require that}
\begin{align}
  & \beta(0) = \beta''(0) = 0
  \aand
  \beta'''(s) \geq 0 
  \quad \hbox{for every $s\in\erre$}
  \label{hpbeta}
  \\
  \separa
  & \lim_{|s|\to+\infty} \frac {\beta'(s)} {|s|} 
  = + \infty
  \label{hpbetaprimo}
  \\[2mm]
  \separa
  & |\beta''(s)|
  \leq \Cbeta ( |\beta'(s)| + 1 )
  \quad \hbox{for some $\Cbeta>0$ and every $s\in\erre$}.
  \label{hpbetasecondo}
\end{align}
\Accorpa\HPstruttura hpconstants hpbetasecondo
Finally
\Beq
  \tau \in (0,1)
  \label{hptau}
\Eeq
is the relaxation parameter.

\Brem
\label{Remstructure}
Let us observe that our structural assumptions imply that \pier{$\beta$ and $\beta''$ are (maximal) monotone,
$\beta'$ and $\beta'''$ are nonnegative} and that the estimates
\begin{align}
  & 0 \leq \beta'(r)
  = \beta'(0) + \int_0^r \beta''(s) \, ds
  \non
  \\
  & \leq \beta'(0) + \Cbeta \pier{ \biggl| \int_0^r (\beta' (s) + 1) \, ds \biggr| }
 = \Cbeta \pier{ |\beta(r)+r| }+ \Cbeta'
 \label{stimabetaprimo}
 \\
 \noalign{\noindent and}
 & |\Beta(r)|
 \leq |\Beta(0)| + \pier{ \int_0^r \beta(s) \, ds }
 \leq \Cbeta'' + |r| \, |\beta(r)|
 \label{stimaBeta}
\end{align}
hold true for every $r\in\erre$ with an obvious meaning of $\Cbeta'$ and~$\Cbeta''$.
Moreover, we remark that our assumptions are satisfied by a wide class of smooth potentials
with either polynomial or exponential growth.
In particular, they are met in the case of the classical potential~\eqref{regpot}.
\Erem

As for the forcing term and the initial data, we assume that
\Beq
  g \in \L2H , \quad
  \phiz \in W 
  \aand
  \rhoz \in H \,.
  \label{hpdati}
\Eeq

At this point, we are ready to make the notion of solution to the relaxed problem precise and to state our existence result.
The variational equation we are going to write can be formally obtained by standard integration by parts
on account of the homogeneous Neumann boundary conditions
starting from equations \eqref{Iprimatau}, \eqref{Iseconda} and \eqref{Iterza} tested by the arbitrary test function~$v$.
Indeed, the solution we are looking for is a triplet $\soluz$ with the regularity
\begin{align}
  & \phi \in \H2\Zp \cap \W{1,\infty}\Vp \cap \L\infty W
  \label{regphi}
  \\
  & \mu \in \L\infty\Wp
  \label{regmu}
  \\
  & w \in \L\infty H
  \label{regw}
\end{align}
\Accorpa\Regsoluz regphi regw
that solves the equations
\begin{align}
  & \tau \< \dt^2\phi , v >
  + \< \dt\phi , v >
  + \sigma \iO \phi v
  - \< \mu , \Delta v >
  = \iO g v 
  \non
  \\
  & \quad \hbox{for every $v\in Z$ and \aet}
  \label{prima}
  \\
  \separa
  & - \iO w \, \Delta v
  + \iO \beta'(\phi) w v
  + (\nu-\lambda) \iO w v
  = \< \mu , v >
  \non
  \\
  & \quad \hbox{for every $v\in W$ and \aet}
  \label{seconda}
  \\
  \separa
  & - \Delta\phi 
  + \beta(\phi)
  - \lambda\phi 
  = w
  \quad \aeQ
  \label{terza}
\end{align}
and satisfies the initial conditions
\Beq
  \phi(0) = \phiz
  \aand
  \dt\phi(0) = \rhoz \,.
  \label{cauchy}
\Eeq
\Accorpa\Pbl prima cauchy
Notice that $\phi$ takes values in~$W$ (see \eqref{regphi})
so that it satisfies the homogeneous Neumann boundary condition as well (see \eqref{defspazi}).
On the contrary, due to the very low level of regularity we are requiring,
equations \eqref{seconda} and \eqref{terza} \pier{represent only a very weak formulation of their 
counterparts}~\eqref{Iprimatau} and \eqref{Iseconda} \pier{with} the boundary conditions for $\mu$ and~$w$. 
\pier{Consequently, this lack of regularity leaves the uniqueness of the solution as an open question.
We now present our existence result.}
\Bthm
\label{Existence}
Assume \HPstruttura\ on the structure of the system and \eqref{hpdati} on the data.
Then, for every $\tau\in(0,1)$, there exists at least one triplet $\soluz$ 
with the regularity specified by \Regsoluz\ that solves Problem \Pbl\
and satisfies the estimates
\begin{align}
  & \norma\phi_{\H1\Vp\cap\L\infty W}
  + \norma\mu_{\L\infty\Wp}
  + \norma w_{\L\infty H}
  \leq C
  \label{stab}
  \\
  & \tau \norma{\dt^2\phi}_{\L2\Zp}
  + \tau^{1/2} \norma{\dt\phi}_{\L\infty\Vp}
  \leq C
  \label{stabbis}
\end{align}
\Accorpa\Stability stab stabbis
with a constant $C$ that depends only on $\Omega$, $T$, the structure of the system and the norms of the data corresponding to \eqref{hpdati}.
In particular, $C$~does not depend on~$\tau$.
\Ethm

\Brem
\label{Continuous}
We notice that $\H1\Vp\cap\L\infty W$ is continuously embedded in $\CQ0$
(see, e.g., \cite[Sect.~8, Cor.~4]{Simon}).
Therefore, $\phi$ is continuous on $\overline Q$ as well as $\beta(\phi)$ and $\beta'(\phi)$
and the maximum norm of these functions are bounded by a constant 
that has the same dependence as the above $C$ has.
Moreover, by a time integration in \eqref{prima},
one could prove one more estimate, namely,
\Beq
  \norma{1*\mu}_{\pier{\L\infty V}} \leq C'
  \non
\Eeq
with a constant $C'$ similar to~$C$, where $1*\mu$ is the function defined~by
\Beq
  (1*\mu)(t) := \iot \mu(s) \, ds
  \quad \hbox{for every $t\in[0,T]$}.
  \non
\Eeq
However, we do not give a detailed proof of this since we do not use it in the sequel.
\Erem

\pier{Now, we present our next result concerning the limit of solutions to problem \Pbl\ as $\tau$ tends to zero.
To this end, we consider the problem with $\tau=0$ as formulated and solved in~\cite{CGSS6}.
Specifically, it involves finding a triplet} $\soluz$ with the regularity 
\begin{align}
  & \phi \in \H1\Vp \cap \L\infty W
  \label{regphiz}
  \\
  & \mu \in \L2V
  \label{regmuz}
  \\
  & w \in \L2W
  \label{regwz}
\end{align}
\Accorpa\Regsoluz regphi regw
that solves the equations
\begin{align}
  & \< \dt\phi , v >
  + \iO \nabla\mu \cdot \nabla v
  + \sigma \iO \phi v
  = \iO g v
  \quad \hbox{for every $v\in V$ and \aet}
  \label{primaz}
  \\
  \separa
  & - \Delta w + \beta'(\phi) w + (\nu-\lambda) w
  = \mu
  \quad \aeQ
  \label{secondaz}
  \\
  & - \Delta\phi + \beta(\phi) - \lambda\phi
  = w
  \quad \aeQ
  \label{terzaz}
\end{align}
and satisfies the initial condition
\Beq
  \phi(0) = \phiz \,.
  \label{cauchyz}
\Eeq
\Accorpa\Pblz primaz cauchyz
Under our assumptions on the structure of the system and the \pier{given data, this problem admits a unique solution, primarily due to \cite[Thm.~2.2]{CGSS6}. 
In that work, the forcing term~$g$ (denoted as $u$ therein)
is assumed to be bounded, as the authors focus on a control problem. 
However, a closer examination of the proof reveals that our assumption on $g$ is sufficient to ensure well-posedness.}

We now present our result.
\Bthm
\label{Asymptotics}
Under the assumption of Theorem~\ref{Existence},
let $\graffe\soluzt$ be any family of solutions to problem \Pbl that satisfies the estimates \Stability.
Then, this family converges to the unique solution \pier{$(\phi,\mu, w)$} to Problem \Pblz\ as $\tau$ tends to zero
in the following sense
\begin{align}
  & \phit \to \phi
  && \quad \hbox{weakly star in $\H1\Vp\cap\L\infty W$}
  \non
  \\
  &&& \hbox{\phantom{\qquad} and uniformly in $\overline Q$}
  \label{convphit}
  \\
  \separa
  & \mut \to \mu
  && \quad \hbox{weakly star in $\L\infty\Wp$}
  \label{convmut}
  \\
  & \wt \to w
  && \quad \hbox{weakly star in $\L\infty H$} \,.
  \label{convwt}
\end{align}
\Ethm

\pier{We can also prove an error estimate in weaker norms for the difference between $\soluzt$ and the limit $(\phi,\mu, w)$ in terms of the parameter $\tau$.}

\pier{\Bthm
\label{ErrEst}
Under the same assumption \gianni{and notation} as in Theorem~\ref{Asymptotics}, 
\gianni{the following estimate}
\begin{align}
  & \norma{\phi_\tau-\phi}_{\C0{V^*}\cap\L2{W}}
  + \norma{\mu_\tau-\mu}_{\L2{W^*}}
  \non
  \\
  & \quad {}
  + \norma{w_\tau-w}_{\L2H}
  \leq C \, \tau^{1/2}
  \label{error}
\end{align}
\gianni{holds true for every $\tau\in(0,1)$ with a constant $C$ that}
depends only on $\Omega$, $T$, the structure of the system and the norms of the data corresponding to \eqref{hpdati}.
\Ethm}%

The rest of the section is devoted to the collection of some useful tools.
First of all, besides the \Holder\ and Cauchy--Schwarz inequalities, we often account for the Young and Poincar\'e inequalities
\begin{align}
  & ab \leq \delta a^2 + \frac 1{4\delta} \, b^2
  \quad \hbox{for every $a,b\in\erre$ and $\delta>0$}
  \label{young}
  \\[2mm]
  & \normaV v
  \leq \CO \, \bigl( \norma{\nabla v} + |\vbar| \bigr)
  \quad \hbox{for every $v\in V$}
  \label{poincare}
\end{align}
where $\vbar$ denotes the mean value of~$v$
and the constant $\CO$~depends only on~$\Omega$.
\pier{We generally define the generalized mean value $\vbar$ of any element} $v\in\Zp$ 
by setting
\Beq
  \vbar := \frac 1{|\Omega|} \, \< v , 1 > 
  \label{defmean}
\Eeq 
where $1$ stands for the constant function that takes the value $1$ everywhere in~$\Omega$.
In the sequel, for any $s\in\erre$, the same symbol $s$ denotes the corresponding constant functions in $\Omega$ and~$Q$.
Notice that the above definition \eqref{defmean} is meaningful, since $1$ actually belongs to~$Z$,
and that $\vbar$ reduces to the usual mean value when $v\in H$.
The same notation $\vbar$ is also used if $v$ is a time-dependent function.

Next, we recall an important tool which is often used in connection with the Cahn--Hilliard \pier{systems.}
Namely, for a given $\zeta\in\Vp$, we consider the generalized Neumann problem of finding
\begin{align}
  z \in V
  \quad \hbox{such that} \quad
  \iO \nabla z \cdot \nabla v
  = \< \zeta , v >
  \quad \hbox{for every $v\in V$}.
  \label{neumann}
\end{align}
Since $\Omega$ is connected and smooth, it is well known that the above problem admits solutions $z$ if and only if $\zeta$ has zero mean value.
Hence, we can introduce the following solution operator $\calN$ by setting
\begin{align}
  & \calN: \dom(\calN) := \graffe{\zeta\in\Vp:\ \zetabar=0} \to \graffe{z\in V:\ \zbar=0},
  \quad 
  {\cal N}: \zeta \mapsto z,
  \label{defN}
\end{align}
where $z$  is the unique solution to \eqref{neumann} that also satisfies $\zbar=0$.
It turns out that $\calN$ is an isomorphism between the above spaces and that the formula
\Beq
  \normaVp\zeta^2 := \norma{\nabla\calN(\zeta-\zetabar)}^2 + |\zetabar|^2
  \quad \hbox{for every $\zeta\in\Vp$}
  \label{normaVp}
\Eeq
defines a Hilbert norm in $\Vp$ that is equivalent to the standard dual norm of~$\Vp$.
From the above properties, one can obtain the following identities:
\begin{align}
  & \iO \nabla\calN\zeta \cdot \nabla v
  = \< \zeta , v >
  \quad  \hbox{for every $\zeta\in\dom(\calN)$ and $v\in V$}
  \label{dadefN}
  \\
  \separa
  & \calN\zeta \in W
  \aand
  -\Delta\calN\zeta 
  = \zeta
  \quad \hbox{for every $\zeta\in\dom(\calN)\cap H$}
  \label{laplaceN}
  \\
  \separa
  & \< \zeta , \calN\zeta >
  = \iO |\nabla\calN\zeta|^2
  = \normaVp\zeta^2
  \quad \hbox{for every $\zeta\in\dom(\calN)$.}
  \label{danormaVp}
\end{align}
\pier{We also point out that}
\Beq
  \< \dt\zeta(t) , \calN\zeta(t) >
  = \frac 12 \, \ddt \, \normaVp{\zeta(t)}^2
  \quad \aat
  \label{propN} 
\Eeq
which holds true for every $\zeta\in\H1\Vp$ satisfying $\zetabar=0$ \aet.
\Accorpa\PropN dadefN propN

\pier{Moreover, we introduce another useful inequality for the sequel. Recalling the compact embeddings $W\emb V$ and $W\emb\CQ0$,  and using elliptic regularity theory, one can easily derive the inequality}
\Beq 
  \normaV v
  \leq \delta \, \norma{\Delta v} + C_{\Omega,\delta} \, \normaVp v
  \quad \hbox{for every $v\in W$ and $\delta>0$}
  \label{compact}
\Eeq
where $\normaVp\cpto$ is the norm \pier{in $\Vp$ defined in \eqref{normaVp}}
and $C_{\Omega,\delta}$ depends on $\Omega$ and~$\delta$.

\pier{We conclude this section by introducing the following convention: the lowercase symbol 
 $c$ denotes a generic constant
that depends only on $\Omega$, $T$, the structure of the system, and the data given in~\eqref{hpdati}. 
In particular, $c$~is independent of $\tau$ and the parameter $n$ introduced in the next section.
Clearly, since $c$ stands for a generic constant, its value may change from line to line and even within the same line.
Furthermore, whenever a positive constant $\delta$ appears in a computation, we use the notation~$\cdelta$ 
instead of a general~$c$, indicating that these constants also depend on~$\delta$. This convention 
does not apply to specific constants referenced explicitly, such as in~\eqref{stab}, 
where a capital letter is used instead of the generic~$c$.}


\section{Existence}
\label{EXISTENCE}
\setcounter{equation}{0}

In this section, we prove Theorem~\ref{Existence}.
To this end, similarly to~\cite{CGSS6},
we give an alternative formulation of \Pbl\ and we solve the new problem 
by introducing a Faedo--Galerkin scheme. 
Here is the new formulation:

\smallskip\noindent
{\bf Equivalent problem:\quad\sl Find a pair $\soluzbis$ with the regularity \accorpa{regphi}{regmu}
that solves \eqref{prima} and the variational equation
\begin{align}
  & - \iO \bigl(
    - \Delta\phi + \beta(\phi) - \lambda\phi  
  \bigr) \Delta v
  + \iO \beta'(\phi) \bigl(
    - \Delta\phi + \beta(\phi) - \lambda\phi  
  \bigr) v
  \non
  \\
  & \quad {}
  + (\nu-\lambda) \iO \bigl(
    - \Delta\phi + \beta(\phi) - \lambda\phi  
  \bigr) v
  \non
  \\
  & = \< \mu , v >
  \quad \hbox{for every $v\in W$ and \aet}
  \label{quarta}
\end{align}
and satisfies the initial conditions \eqref{cauchy}}.

\smallskip

Clearly, if $\soluz$ is a solution to \Pbl, then \eqref{quarta}~follows for $\soluzbis$
by just replacing $w$ given by \eqref{terza} in \eqref{seconda}.
Conversely, if $\soluzbis$ solves the new problem,
by taking \eqref{terza} as a definition of~$w$,
\eqref{regw}~and \eqref{seconda} immediately follow from \eqref{regphi} and~\eqref{quarta}
(see also \cite[Prop~3.1]{CGSS6} for more details in a similar situation).

As just said, we solve the new problem by a Faedo--Galerkin scheme.
To this end, we introduce the sequence $\graffe{\lambdaj}_{j\geq1}$ of the eigenvalues 
and an orthonormal system $\graffe{\ej}_{j\geq1}$ of corresponding eigenfunctions
of the Neumann problem for the Laplace equation,~i.e.
\begin{align}
  & 0 = \lambda_1 < \lambda_2 \leq \lambda_3 \leq \dots
  \aand \lim_{j\to\infty} \lambdaj = + \infty
  \label{eigenvalues}
  \\
  & \ej\in V
  \aand
  \iO \nabla\ej \cdot \nabla v
  = \lambda_j \iO \ej v
  \quad \hbox{for every $v\in V$ and $j=1,2,\dots$}
  \label{eigenfunctions}
  \\
  & \iO \ei \ej = \delta_{ij}
  \quad \hbox{for $i,j=1,2,\dots$}
  \aand 
  \hbox{$\graffe{e_j}_{j\geq1}$ is a complete system in $H$}
  \label{orthogonality}
\end{align}
\Accorpa\EigenPbl eigenvalues orthogonality
where $\delta_{ij}$ is the Kronecker symbol.
We~set
\Beq
  \Vn := \Span \graffe{e_1,\dots,e_n},
  \quad \hbox{for $n=1,2,\dots$}
  \label{defVn}
\Eeq
and observe that each $\Vn$ is included in $Z$ and that the union of these spaces is dense in both $V$ and~$H$.
We also notice that $V_1=\Span \graffe{e_1}$ consists of the space of constant functions
since $\Omega$ is connected. 

\Brem
\label{Projecion}
Let us make some observations regarding the orthogonal projection operator $\Pn:H\to\Vn$ 
(orthogonality being understood with respect to the standard inner product of~$H$).
Let $Y$ be any of the spaces $H$, $V$, $W$ and~$Z$.
Then, for every $v\in Y$, we have~that 
\Beq
  \norma{\Pn v}_Y \leq \CO \, \norma v_Y
  \aand
  \Pn v \to v
  \quad \hbox{strongly in $Y$}
  \label{convPn}
\Eeq
where the constant~$\CO$ depends only on~$\Omega$.
Notice that the convergence in \eqref{convPn} actually is strong by uniform convexity
if we can take $\CO=1$ by carefully choosing the norm of~$Y$.
Clearly, this is the case if $Y=H$ with its natural norm.
The same is true if $Y=V$ with its standard Hilbert norm
(see \cite[Rem~3.3]{CGSS1} for details).
It follows from the former observation that we can tale $\CO=1$
if either $Y=W$ or $Y=Z$ as well, 
provided that their norms are the graph norms defined~by
\Beq
  \normaW v^2 = \norma v^2 + \norma{\Delta v}^2
  \aand
  \normaZ v^2 = \normaW v^2 + \normaW{\Delta v}^2
  \non
\Eeq
for $v\in W$ and $v\in Z$, respectively.
This immediately follows from the formula $\Delta\Pn v=\Pn\Delta v$ for $v\in W$,
which we prove at once.

Let $v\in W$.
Then $\Delta\Pn v\in\Vn$ since $\Pn v\in\Vn$.
On the other hand,
\Beq
  (\Delta\Pn v,z)
  = (\Pn v,\Delta z)
  = (v,\Delta z)
  = (\Delta v,z)
  \quad \hbox{for every $z\in W$}
  \non
\Eeq
since $\Pn v\in W$.
As $W$ is dense in~$H$, it follows that
$(\Delta\Pn v,z)=(\Delta v,z)$ for every $z\in H$,
that is, $\Delta\Pn v=\Pn\Delta v$.

Next, if $v\in\L2Y$ and $\vn$ is defined by
$\vn(t):=\Pn(v(t))$ \aat\ (and we simply write $\vn=\Pn v$ in the sequel),
then
\Beq
  \norma\vn_{\L2Y} \leq \CO \, \norma v_{\L2Y},
  \aand
  \vn \to v 
  \quad \hbox{strongly in $\L2Y$}
  \label{convPnbis}
\Eeq
where one can take $\CO=1$ if the norm in $Y$ is well chosen as said before.
\Erem

At this point, we are ready to introduce the Faedo--Galerkin scheme mentioned above.
It is the following

\smallskip\noindent
{\bf Problem \Pbln:\quad\sl Look for a pair 
\Beq
  \soluzbis \in \H2\Vn \times \L2\Vn
  \label{regsoluzn}
\Eeq
that solves the variational equations \eqref{prima} and \eqref{quarta} just for test functions $v\in\Vn$
and satisfies the initial conditions 
\Beq
  \phi(0) = \Pn\phiz
  \aand
  \phi'(0) = \Pn\rhoz \,.
  \label{cauchyn}
\Eeq
}

\step
Solution to the discrete problem

The discrete problem just introduced is equivalent to a second order Cauchy problem for a system of ordinary differential equations,
as we briefly see.
We expand the unknown $\phi$ and~$\mu$ in terms of the eigenfunction $\ej$ as follows
\Beq
  \phi(t) = \somma j1n \phij(t) \ej
  \aand
  \mu(t) = \somma j1n \muj(t) \ej \, .
  \non
\Eeq
Then, the new unknowns are the functions 
\Beq
  \bphi := (\phij)_{j=1}^n \in \H2\erren
  \aand
  \bmu := (\muj)_{j=1}^n \in \L2\erren 
  \non
\Eeq
and we can rewrite the variational equations \eqref{prima} and~\eqref{quarta} 
(which have to be satisfied for every~$v\in\Vn$)
in~terms of them.
Since we can just take $v=\ei$ for $i=1,\dots,n$ instead of every $v\in\Vn$, we have for $i=1,\dots,n$
\begin{align}
  & \tau \phii''
  + \phii'
  + \sigma \phii
  + \lambdai \mui 
  = \iO g \ei
  \label{primaij}
  \\
  & \iO \omega \lambdai \ei
  + \iO \gamma \omega \lambdai \ei
  + (\nu-\lambda) \iO \omega \ei
  = \mui
  \label{quartaij}
\end{align} 
where we have set for brevity
\Beq
  \omega := \somma j1n (\lambdaj-\lambda) \phij \ej + \beta \Bigl( \somma j1n \phij\ej \Bigr)
  \aand
  \gamma := \beta' \Bigl( \somma j1n \phij\ej \Bigr) \,.
  \non
\Eeq
Clearly, we can \pier{replace $\mui$ in \eqref{primaij} by using the expression in} \eqref{quartaij}.
Therefore, the system to be solved takes the~form
\Beq
  \tau \bphi'' + \bphi' + \calF(\bphi) = \gg
  \non
\Eeq
where $\calF:\erren\to\erren$ is of class~$C^2$ and $\gg\in\L2\erren$.
Moreover, the initial conditions \eqref{cauchyn} provide prescribed values for $\bphi(0)$ and~$\bphi'(0)$
and we obtain a well posed Cauchy problem.
Therefore, it has a unique maximal solution $\bphi$ that is defined in the interval $[0,T_n)$ for some $T_n\in(0,T]$.
Consequently, the discrete problem under investigation has a unique maximal solution $\soluzn$ defined in~$[0,T_n)$.

\smallskip

However, the estimates we are going to perform show that the discrete solution is bounded, so that $T_n=T$ by maximality.
For that reason, to simplify the notation, we already write $T$ in place of~$T_n$.
We will come back to this point later on. 
Moreover, the estimates we prove are also sufficient to let $n$ tend to infinity 
and construct a solution $\soluzbis$ of the new version of Problem \Pbl,
thus a solution $\soluz$ to \Pbl\ itself, with the desired stability estimates.
In order to simplify the notation, we write $\soluzbis$ in place of $\soluzn$ in performing the estimates
and restore the notation $\soluzn$ only at the end of each procedure.
Moreover, for a given time dependent test function~$v$,
it is understood that the equations we consider are always written at the time $t$ and then tested by $v(t)$ \aat\ 
even though we avoid writing the time, for simplicity.
In the sequel, we use the notation
\Beq
  Q_t := \Omega \times (0,t)
  \quad \hbox{for $t\in(0,T]$} \,.
  \label{defQt}
\Eeq

\step
First a priori estimate

We recall the notation \eqref{defmean} and that all the constant functions belong su~$V_1$, thus to any~$\Vn$.
Hence, we can test \eqref{prima} by $1/|\Omega|$ and obtain
\Beq
  \tau \phibar'' + \phibar' + \sigma\phibar
  = \gbar 
  \quad \aet \,.
  \label{meanprima}
\Eeq 
Moreover, $\phibar(0)=\phizbar$ and $\phibar'(0)=\rhozbar$ for $(\Pn v)_\Omega=\vbar$ for every $v\in\Vn$.
Since \eqref{hpdati} implies that $\gbar\in L^2(0,T)$, by testing \eqref{meanprima} by~$\phibar'$ \pier{and integrating over $(0,t)$, $t\in (0,T]$,} we easily obtain
\Beq
  \tau^{1/2} \, \norma{(\phin)_\Omega'}_{L^\infty(0,T)}
  + \norma{(\phin)_\Omega}_{H^1(0,T)} \leq c \,.
  \label{primastima}
\Eeq

\step
Second a priori estimate

We take the difference between \eqref{prima} and \eqref{meanprima} 
multiplied by $v\in\Vn$ and integrated over~$\Omega$.
By recalling the definitions \eqref{defN} and \eqref{normaVp} 
of the operator $\calN$ and of the norm $\normaVp\cpto$
as well as the properties \PropN\ of~$\calN$,
we test the difference just constructed by $\calN(\dt\overline{\phi} - \overline{\phi}_\Omega')$
and integrate with respect to time.
We observe that this is correct since even $\calN(\dt\overline{\phi} - \overline{\phi}_\Omega')$ takes values in~$\Vn$ 
(see \cite[Rem.~3.3]{CGSS6} for details).
We~obtain for every $t\in[0,T]$
\begin{align}
  & \frac\tau2 \, \normaVp{(\dt\overline{\phi} - \overline{\phi}_\Omega')(t)}^2
  + \iot \normaVp{(\dt\overline{\phi} - \overline{\phi}_\Omega')(s)}^2 \, ds
  + \frac\sigma2 \, \normaVp{(\overline{\phi} - \overline{\phi}_\Omega)(t)}^2
  \non
  \\
  & = \frac\tau2 \, \normaVp{\rhoz-\rhozbar}^2
  + \frac\sigma2 \, \normaVp{\phiz-\phizbar}^2
  + \intQt (g-\gbar) \, \calN(\dt\overline{\phi} - \overline{\phi}_\Omega') 
  \non
  \\
  & \quad {}
  - \intQt \mu \, (\dt\overline{\phi} - \overline{\phi}_\Omega') \,.
  \label{testprima}
\end{align}
Next, we consider the energy $\calE(\phi(t))$, i.e., \eqref{Ienergy} with $v=\phi(t)$.
However, we avoid writing the time $t$ for a while.
By accounting for \eqref{calcolof}, we~have
\begin{align}
  & \calE(\phi)
  = \frac12  \iO \bigl( -\Delta\phi + \beta(\phi) - \lambda\phi \bigr)^2
  + \nu \iO \Bigl( \frac 12 \, |\nabla\phi|^2 + \Beta(\phi) - \frac\lambda2 \, \phi^2 \Bigr).
  \non
\end{align}
From one side, we rewrite this by splitting the first integral into two parts
and making explicit calculations in the second one.
On the other hand, we compute the time derivative of the energy.
We have
\begin{align}
  & \calE(\phi)
  = \frac14  \iO \bigl( -\Delta\phi + \beta(\phi) - \lambda\phi \bigr)^2
  + \frac14 \iO \bigl(
    |\Delta\phi|^2 + |\beta(\phi)|^2 + \lambda^2 |\phi|^2
  \bigr)
  \non
  \\
  & \quad {}
  + \frac 12 \iO \bigl(
    \beta'(\phi) |\nabla\phi|^2 - \lambda |\nabla\phi|^2 - \lambda \beta(\phi) \phi
  \bigr)
  + \nu \iO \Bigl( \frac 12 \, |\nabla\phi|^2 + \Beta(\phi) - \frac\lambda2 \, \phi^2 \Bigr)
  \label{energy}
  \\
  \separa
  \noalign{\noindent as well as \smallskip}
  & \ddt \, \calE(\phi)
  = \iO \bigl(
    -\Delta\phi + \beta(\phi) - \lambda\phi
  \bigr)
  \bigl(
    -\Delta\dt\phi + \beta'(\phi)\dt\phi - \lambda\dt\phi
  \bigr)  
  \non
  \\
  & \quad {}
  + \nu \iO \bigl( 
    \nabla\phi \cdot \nabla\dt\phi + \beta(\phi) \dt\phi - \lambda \phi \dt\phi
  \bigr).
  \label{dtenergy}
\end{align}
At this point, we test \eqref{quarta} by $\dt\overline{\phi} - \overline{\phi}_\Omega'$ and rearrange a little.
We~have
\begin{align}
  & \iO \bigl(
    -\Delta\phi + \beta(\phi) - \lambda\phi
  \bigr) (-\Delta\dt\phi)
  + \iO \beta'(\phi) \bigl(
    - \Delta\phi + \beta(\phi) - \lambda\phi
  \bigr) \dt\phi
  \non
  \\
  & \quad {}
  - \lambda \iO \bigl(
    -\Delta\phi + \beta(\phi) - \lambda\phi
  \bigr) \dt\phi
  + \nu \iO \bigl(
    - \Delta\phi + \beta(\phi) - \lambda\phi
  \bigr) \dt\phi
  \non
  \\
  \separa
  & = \iO \beta'(\phi) \bigl(
    - \Delta\phi + \beta(\phi) - \lambda\phi
  \bigr) \phibar'
  + (\nu-\lambda) \iO \bigl(
    - \Delta\phi + \beta(\phi) - \lambda\phi
  \bigr) \phibar'
  \non
  \\
  & \quad {}
  + \iO \mu (\dt\overline{\phi} - \overline{\phi}_\Omega')
  \non
\end{align}
and we see that the \lhs\ of this equality coincides with the \rhs\ of \eqref{dtenergy}.
Hence, if we integrate over $(0,t)$, we obtain
\Beq
  \calE(\phi(t))
  = \calE(\phiz)
  + \iot R(s) \, ds
  + \intQt \mu \, (\dt\overline{\phi} - \overline{\phi}_\Omega') 
  \label{testquarta}
\Eeq
where we have set for brevity
\Beq
  R :=
  \iO (\beta'(\phi)+\nu-\lambda) \bigl(
    - \Delta\phi + \beta(\phi) - \lambda\phi
  \bigr) \phibar' \,.
  \label{defR}
\Eeq
At this point, we add \eqref{testprima} and \eqref{testquarta} to each other
and notice an obvious cancellation.
By using the form \eqref{energy} of the energy,
we conclude that
\begin{align}
  & \frac\tau2 \, \normaVp{(\dt\overline{\phi} - \overline{\phi}_\Omega')(t)}^2
  + \iot \normaVp{(\dt\overline{\phi} - \overline{\phi}_\Omega')(s)}^2 \, ds
  + \frac\sigma2 \, \normaVp{(\overline{\phi} - \overline{\phi}_\Omega)(t)}^2
  \non
  \\
  & \quad {}
  + \frac14  \iO \bigl|\bigl( -\Delta\phi + \beta(\phi) - \lambda\phi \bigr)(t)\bigr|^2
  + \frac14 \iO \bigl(
    |\Delta\phi|^2 + |\beta(\phi)|^2 + \lambda^2 |\phi|^2
  \bigr)(t)
  \non
  \\
  & \quad {}
  + \frac 12 \iO \bigl(
    \beta'(\phi) |\nabla\phi|^2 - \lambda |\nabla\phi|^2 - \lambda \beta(\phi) \phi
  \bigr)(t)
  + \nu \iO \Bigl( \frac 12 \, |\nabla\phi|^2 + \Beta(\phi) - \frac\lambda2 \, \phi^2 \Bigr)(t)
  \non
  \\ 
  \separa
  & = \frac\tau2 \, \normaVp{\rhoz-\rhozbar}^2
  + \frac\sigma2 \, \normaVp{\phiz-\phizbar}^2
  + \intQt (g-\gbar) \, \calN(\dt\overline{\phi} - \overline{\phi}_\Omega') 
  \non
  \\
  & \quad {}
  + \calE(\phiz)
  + \iot R(s) \, ds \,.
  \label{persecondastima}
\end{align}
From this identity we derive our basic estimate.
We keep just the nonnegative terms on the \lhs,
move the others to the \rhs\ and estimate 
both them and those that already form the \rhs\ of \eqref{persecondastima}.
What remains on the \lhs\ is the following
\begin{align}
  & \frac\tau2 \, \normaVp{(\dt\overline{\phi} - \overline{\phi}_\Omega')(t)}^2
  + \iot \normaVp{(\dt\overline{\phi} - \overline{\phi}_\Omega')(s)}^2 \, ds
  + \frac\sigma2 \, \normaVp{(\overline{\phi} - \overline{\phi}_\Omega)(t)}^2
  \non
  \\
  & \quad {}
  + \frac14  \iO \bigl|\bigl( -\Delta\phi + \beta(\phi) - \lambda\phi \bigr)(t)\bigr|^2
  + \frac14 \iO \bigl(
    |\Delta\phi|^2 + |\beta(\phi)|^2 + \lambda^2 |\phi|^2
  \bigr)(t)
  \non
  \\
  & \quad {}
  + \frac 12 \iO \bigl( \beta'(\phi) |\nabla\phi|^2 \bigr)(t) \,.
  \label{lhssecondastima}
\end{align}
Now, we estimate the terms on the new \rhs, separately.
In our computation, we repeatedly owe to the Cauchy--Schwarz and Young inequalities (see~\eqref{young}).
Moreover, $\delta$~is a positive parameter whose value is chosen later on.
In order to treat the first term we consider, 
we account for \eqref{compact} and denote by $(\cpto,\cpto)_*$ the inner product in $\Vp$ 
associated with the norm $\normaVp\cpto$ given by \eqref{normaVp}.
We~have
\begin{align}
  & \frac {\lambda-\nu}2 \iO |\nabla\phi(t)|^2
  \leq \delta \iO |\Delta\phi(t)|^2
  + \cdelta \, \normaVp{\phi(t)-\phibar(t)}^2
  \non
  \\
  \separa
  & \leq \delta \iO |\Delta\phi(t)|^2
  + \cdelta \, \normaVp{\phiz-\phizbar}^2
  + \cdelta \iot \bigl(
    \phi(s)-\phibar(s) , \dt\phi(s)-\phibar'(s)
  \bigr)_* \, ds
  \non
  \\  
  & \leq \delta \iO |\Delta\phi(t)|^2
  + \cdelta
  + \delta \iot \normaVp{\dt\phi(s)-\phibar'(s)}^2 \, ds
  + \cdelta \iot \normaVp{\phi(s)-\phibar(s)}^2 \, ds \,.
  \non
\end{align}
Next, we have
\begin{align}
  &  \frac \lambda2 \iO \beta(\phi(t)) \phi(t)
  \leq \delta \iO |\beta(\phi(t)|^2 
  + \cdelta \iO |\phi(t)|^2
  \non
\end{align}
and the last term can be dealt with as above, 
since the inequality \eqref{compact} we have used there involves the full $V$-norm.
The same integral occurs in the next estimate.
\pier{With the help of}~\eqref{stimaBeta} we have indeed
\begin{align}
  & - \nu \iO \Bigl( \frac 12 \, \Beta(\phi) - \frac\lambda2 \, \phi^2 \Bigr)(t)
  \leq \delta \iO |\beta(\phi(t))|^2
  + \cdelta \iO |\phi(t)|^2
  + \cdelta \,.
  \non
\end{align}
Since the terms involving the initial data are bounded thanks to our assumptions~\eqref{hpdati},
we consider the integral related to the forcing term. 
We recall the properties of $\calN$ already mentioned and~have
\begin{align}
  & \intQt (g-\gbar) \, \calN(\dt\overline{\phi} - \overline{\phi}_\Omega') 
  \leq c \iot \normaVp{g(s)} \, \normaV{\calN(\dt\phi(s)-\phibar'(s))} \, ds
  \non
  \\
  & \leq c \iot \norma{g(s)} \, \normaVp{\dt\phi(s)-\phibar'(s)} \, ds
  \leq \delta \iot \normaVp{\dt\phi(s)-\phibar'(s)}^2\, ds
  + \cdelta \,.
  \non
\end{align}
Finally, the term involving the function $R$\pier{, defined in \eqref{defR}, can be handled using \eqref{stimabetaprimo} as follows:}
\begin{align}
  & \iot R(s) \, ds 
  \leq c \intQt \bigl( 1 + |\phi| + |\beta(\phi)| \bigr) \, |{-\Delta\phi}+\beta(\phi)-\lambda\phi| \, |\phibar'|
  \non
  \\
  & \leq c \iot |\phibar'(s)| \bigl(
    1 + \norma{\phi(s)}^2 + \norma{\beta(\phi(s))}^2
    + \norma{(-\Delta\phi+\beta(\phi)-\lambda\phi)(s)}^2
  \bigr) \, ds .
  \non
\end{align}
\pier{We observe that the last integral above} can be treated by the Gronwall lemma 
since the coefficient $\phibar'$ \pier{is bounded}
in $L^2(0,T)$, thus in $L^1(0,T)$, \pier{due to}~\eqref{primastima}.
By comparing all these inequalities with the \lhs\ given by \eqref{lhssecondastima}, 
choosing $\delta$ small enough and applying the Gronwall lemma,
we obtain an estimate involving~$\phibar'$.
However, this has been already estimated by~\eqref{primastima}.
Therefore, we conclude~that
\begin{align}
  & \tau^{1/2} \, \norma{\dt\phin}_{\L\infty\Vp}
  + \norma\phin_{\H1\Vp\cap\L\infty H}
  \non
  \\
  & \quad {}
  + \norma{\Delta\phin}_{\L\infty H}
  + \norma{\beta(\phin)}_{\L\infty H}
  \leq c \,.
  \label{secondastima}
\end{align}
Thus, elliptic regularity gives
\Beq
  \norma\phin_{\L\infty W} \leq c 
  \quad \hbox{whence also} \quad
  \norma\phin_{\CQ0} \leq c 
  \label{dasecondastima}
\Eeq
since we can apply, e.g., \cite[Sect.~8, Cor.~4]{Simon}).
We have written $T$ in place of~$T_n$ even though this is not correct in principle.
However, \eqref{dasecondastima} written with $T_n$ in place of~$T$
ensures that the maximal solution of the system of ordinary differential equations equivalent to Problem~\pbln\
is bounded too, thus global.
Therefore, a~posteriori, all our estimates are correct as they are.
In particular, \eqref{primastima}, \eqref{secondastima} and \eqref{dasecondastima} are correct.

\step
Third a priori estimate

Let $v\in W$.
We test \eqref{quarta} by $\vn:=\Pn v$.
\pier{By virtue of \eqref{convPn} with $Y=W$, we easily deduce~that}
\Beq
  \iO \mu \vn 
  \leq c \, \bigl(
    \normaW\phi + \norma{\beta(\phi)}_\infty + \norma\phi_\infty
  \bigr) \bigl( 
    1 + \norma{\beta'(\phi)}_\infty
  \bigr) \normaW\vn
  \leq c \, \normaW v
  \non
\Eeq
\aet. 
But $(\mu,\vn)=(\mu,v)$ since $\mu$ takes values in~$\Vn$.
Therefore, we conclude that
\Beq
  \norma\mun_{\L\infty\Wp} \leq c \,.
  \label{terzastima}
\Eeq

\step
Fourth a priori estimate

Let $v\in\L2Z$.
We test \eqref{prima} by $\vn:=\Pn v$.
Similarly as in the above calculation we have~that
\Beq
  \tau \iO \dt^2\phi \, v
  = \tau \iO \dt^2\phi \, \vn
  \leq c \bigl(
    \normaVp{\dt\phi} + \norma\phi + \norma g
  \bigr) \normaV\vn
  + c \, \norma\mu_{\Wp} \, \normaZ\vn 
  \non
\Eeq
\aet.
By squaring, integrating over $(0,T)$ and applying the above estimates and \eqref{convPnbis} with $Y=Z$, we conclude that
\Beq
  \pier{\tau} \norma{\dt^2\phi}_{\L2\Zp} \leq c \,.
  \label{quartastima}
\Eeq

\step
Conclusion

Our aim is letting $n$ tend to infinity.
By the above estimates, well-known compactness results
(for strong compactness see, e.g., \cite[Sect.~8, Cor.~4]{Simon})
and the smoothness of~$\beta$,
there exists a pair $\soluzbis$ such that
\begin{align}
  & \phin \to \phi
  && \quad \hbox{weakly star in $\H2\Zp\cap\W{1,\infty}\Vp\cap\L\infty W$}
  \non
  \\
  &&& \hbox{\phantom{\qquad} and uniformly in $\overline Q$}
  \label{convphin}
  \\
  & \mun \to \mu
  && \quad \hbox{weakly star in $\L\infty\Wp$}
  \label{convmun}
  \\
  \separa
  & \beta(\phin) \to \beta(\phi)
  && \quad \hbox{uniformly in $\overline Q$}
  \label{convbetaphin}
  \\
  & \beta'(\phin) \to \beta'(\phi)
  && \quad \hbox{uniformly in $\overline Q$}
  \label{convbetaprimophin}
\end{align}
at least for a subsequence.
We want to prove that $\soluzbis$ is a solution we are looking for.
By combining the convergence of $\phin(0)$ and $\dt\phin(0)$ at least in $\Zp$
with \eqref{convPn}, we see that \eqref{cauchyn} imply~\eqref{cauchy}.
It remains to prove that $\soluzbis$ solves the variational equations \eqref{prima} and~\eqref{quarta}.
To this end, it suffices to show that $\soluzbis$ solves their time integrated versions, namely,
\begin{align}
  & \tau \< \dt^2\phi , v >_{\L2Z}
  + \< \dt\phi , v >_{\L2V}
  + \sigma \intQ \phi v
  - \< \mu , \Delta v >_{\L2W}
  \non
  \\
  & = \intQ g v 
  \quad \hbox{for every $v\in\L2Z$}
  \label{intprima}
  \\
  \separa
  & - \intQ \bigl(
    - \Delta\phi + \beta(\phi) - \lambda\phi  
  \bigr) \Delta v
  + \intQ \beta'(\phi) \bigl(
    - \Delta\phi + \beta(\phi) - \lambda\phi  
  \bigr) v
  \non
  \\
  & \quad {}
  + (\nu-\lambda) \intQ \bigl(
    - \Delta\phi + \beta(\phi) - \lambda\phi  
  \bigr) v
  \non
  \\
  & = \< \mu , v >_{\L2W}
  \quad \hbox{for every $v\in\L2W$} \,.
  \label{intquarta}
\end{align}
To prove this, let us observe that the solution $\soluzn$ to Problem~\pbln\
solves the integrated version of it,~i.e.,
\begin{align}
  & \tau \intQ \dt^2\phin \, \vn 
  + \intQ \dt\phin \, \vn 
  + \sigma \intQ \phin \vn
  - \intQ \mun \, \Delta\vn 
  \non
  \\
  & = \intQ g \vn
  \quad \hbox{for every $\vn\in\L2\Vn$}
  \label{intpriman}
  \\
  \separa
  & - \intQ \bigl(
    - \Delta\phin + \beta(\phin) - \lambda\phin  
  \bigr) \Delta \vn
  + \intQ \beta'(\phin) \bigl(
    - \Delta\phin + \beta(\phin) - \lambda\phin 
  \bigr) \vn
  \non
  \\
  & \quad {}
  + (\nu-\lambda) \intQ \bigl(
    - \Delta\phin + \beta(\phin) - \lambda\phin  
  \bigr) \vn
  \non
  \\
  & = \intQ \mun \, \vn 
  \quad \hbox{for every $\vn\in\L2\Vn$} \,.
  \label{intquartan}
\end{align}
Then, given $v\in\L2Z$, we can choose $\vn=\Pn v$ in \eqref{intpriman}
and given $v\in\L2W$, we can choose $\vn=\Pn v$ in \eqref{intquartan}.
At this point, we let $n$ tend to infinity and derive \eqref{intprima} and \eqref{intquarta}, respectively. 
Indeed, to justify this procedure, it suffices to combine \accorpa{convphin}{convbetaprimophin}
with the strong convergence~\eqref{convPnbis}.
Finally, the solution we have constructed satisfies the estimates
that follow from \eqref{primastima}, \eqref{secondastima} and \accorpa{terzastima}{quartastima} and the weak \pier{lower semicontinuity} of the norms.
In particular, one can take the same values of the constants $c$ that appear in the estimates we have established.
This proves a part of \eqref{stab} and the full \eqref{stabbis}.
Namely, just the estimate of $w$ is missing.
However, by recovering the triplet $\soluz$ as said at the beginning of the present section,
i.e., by taking \eqref{terza} as a definition of~$w$,
we see that $\soluz$ solves Problem \Pbl\ and satisfies the estimate for $w$ as~well since
\Beq
  \norma w
  \leq \norma{\Delta\phi} + \norma{\beta(\phi)} + \lambda \, \norma\phi
 \quad \aet \,.
  \non
\Eeq
This completes the proof of Theorem~\ref{Existence}.


\section{Asymptotics}
\label{ASYMPTOTICS}
\setcounter{equation}{0}

In this section, we prove the \pier{Theorems~\ref{Asymptotics} and~\ref{ErrEst}.}
We prepare two regularity results that have an independent interest.
The second one, which deals with an elliptic problem in a very weak form, should be known.
Nevertheless, we give a short proof for the reader's convenience.

\Blem
\label{Permu}
Let $u\in\L2\Vp$ satisfy
\Beq
  \ioT \< u(t) , v(t) > \, dt
  \leq M \, \norma v_{\L2\Vp}
  \quad \hbox{for every $v\in\L2W$}
  \label{hppermu}
\Eeq
where $M$ is some positive constant.
Then
\Beq
u \in \L2V.
  \label{permu}
\Eeq
\Elem

\Bdim
We set for brevity
\Beq
  \calV := \L2V 
  \aand
  \calW := \L2W
  \non
\Eeq
so that $u\in\calV^*$ and \eqref{hppermu} becomes
\Beq
  \< u , v >_\calV \leq M \norma v_{\calV^*}
  \quad \hbox{for every $v\in\calW$} \,.
  \non
\Eeq
Since all the embeddings in \eqref{embeddings} are dense,
we have that $\calW$ is dense in~$\calV^*$ 
and the map $\calW\ni v\mapsto\<u,v>_{\calV}$ can be continuously extended to the whole of~$\calV^*$.
Thus, there exists $u^{**}\in\calV^{**}$ such that
\Beq
  \< u^{**} , v >_{\calV^*}
  = \< u , v >_{\calV}
  \quad \hbox{for every $v\in\calW$}.
  \non
\Eeq
Since $\calV$ is reflexive, there exists $\hat u\in\calV$ such that
\Beq
  \< u^{**} , v >_{\calV^*}
  = \< v , \hat u >_\calV
  \quad \hbox{for every $v\in\calV^*$}.
  \non
\Eeq
We thus have that
\Beq
  \< u , v >_\calV
  = \< v , \hat u >_\calV
  \quad \hbox{for every $v\in\calW$}.
  \non
\Eeq
On the other hand, it holds that
\Beq
  \< v , \hat u >_\calV
  = \intQ v \hat u
  = \< \hat u , v >_{\calV}
  \quad \hbox{for every $v\in\calV$}.
  \non
\Eeq
Therefore
\Beq
  \< u , v >_\calV
  = \< \hat u, v >_\calV
  \quad \hbox{for every $v\in\calW$}.
  \non
\Eeq
Since both $u$ and $\hat u$ belong to $\calV^*$ and $\calW$ is dense in $\calV$, 
we conclude that $u=\hat u\in\calV$.
\Edim

\Blem
\label{Perw}
Let $h\in H$ and $u\in H$ and assume that
\Beq
  \iO u (-\Delta v + v) = \iO hv
  \quad \hbox{for every $v\in W$}.
  \label{veryweak}
\Eeq
Then
\Beq
  u \in W \,, \quad
  - \Delta u + u = h
  \quad \aeO 
  \aand
  \normaW u \leq \CO \, \norma h
  \label{perw}
\Eeq
with a constant $\CO$ that depends only on~$\Omega$.
\Elem

\Bdim
It is well known that the problem of finding $\xi$ satisfying
\Beq
  \xi \in W
  \aand
  - \Delta\xi + \xi = h
  \quad \aeO 
  \non
\Eeq
has a unique solution and that $\xi$ satisfies the estimate
$\normaW\xi\leq\CO\,\norma h$ with $\CO$ as in the statement.
Therefore, it suffices to prove that $u=\xi$.
To this end, we observe that $\xi$ also satisfies~\eqref{veryweak}.
Hence, the equality $u=\xi$ follows if we prove that 
the solution $u\in H$ to \eqref{veryweak} is unique.
By linearity, we assume that \eqref{veryweak} holds with~$h=0$
and consider the unique funcion $v$ that satisfies
\Beq
  v \in W
  \aand
  - \Delta v + v = u
  \quad \aeO \,.
  \non
\Eeq
Then, \eqref{veryweak} with~$h=0$ yields
\Beq
  \iO |u|^2
  = \iO u (-\Delta v + v)
  = 0
  \non
\Eeq
whence $u=0$.
\Edim

\pier{\subsection{Proof of Theorem~\ref{Asymptotics}}}%

\pier{We start the proof by picking} a family $\graffe\soluzt$ as in the statement.
Thanks to estimates \Stability\ and well known weak and weak star compactness results,
there exists a triplet $\soluz$ satisfying \accorpa{convphit}{convwt}
as $\tau$ tends to zero.
More precisely, this holds at least for a subsequence, in principle.
However, once it is proved that $\soluz$ satisfies Problem \Pblz,
the whole family $\graffe\soluzt$ converges to $\soluz$ 
since the solution to this problem is unique.
Furthermore, invoking strong compactness as in, e.g., \cite[Sect.~8, Cor.~4]{Simon}
and the smoothness of~$\beta$,
we infer~that
\Beq
  \phit \to \phi \,, \quad
  \beta(\phit) \to \beta(\phi)
  \aand
  \beta'(\phit) \to \beta'(\phi)
  \quad \hbox{uniformly in $\overline Q$}.
  \label{strongphit}
\Eeq
Finally, it is clear that \eqref{convphit} and \eqref{stabbis} imply
\Beq
  \tau \dt^2 \phit \to 0
  \quad \hbox{weakly in $\L2\Zp$} \,.
  \non
\Eeq
Therefore, it is easy to see that $\soluz$ solves Problem~\Pblz\ in a weak form.
Namely, the regularity requirements \accorpa{regmuz}{regwz} of $\mu$ and $w$ are missing 
and just \eqref{terzaz} and \eqref{cauchyz} are satisfied as they are,
while \eqref{primaz} and \eqref{secondaz} are replaced here~by 
\begin{align}
  & \< \dt\phi , v >_{\L2V}
  + \sigma \intQ \phi v
  - \< \mu , \Delta v >_{\L2W}
  \non
  \\
  & = \intQ g v 
  \quad \hbox{for every $v\in\L2Z$}
  \label{intprimaz}
  \\
  \separa
  \noalign{\noindent and}
  & - \intQ w \, \Delta v
  + \intQ \beta'(\phi) \, w \, v
  + (\nu-\lambda) \intQ w v
  \non
  \\
  & = \< \mu , v >_{\L2W}
  \quad \hbox{for every $v\in\L2W$} 
  \label{intsecondaz}
\end{align}
respectively.
However, the above variational equations are equivalent to the following ones
\begin{align}
  & \< \dt\phi , v >
  + \sigma \iO \phi v
  - \< \mu , \Delta v >
  \non
  \\
  & = \iO g v 
  \quad \hbox{for every $v\in Z$ and \aet}
  \label{primazbis}
  \\
  \separa
  \noalign{\noindent and}
  & - \iO w \, \Delta v
  + \iO \beta'(\phi) \, w \, v
  + (\nu-\lambda) \iO w v
  \non
  \\
  & = \< \mu , v >
  \quad \hbox{for every $v\in W$ and \aet} 
  \label{secondazbis}
\end{align}
and we want to recover \eqref{primaz} and \eqref{secondaz}
by proving the regularity requirements \eqref{regmuz} and \eqref{regwz} for $\mu$ and~$w$.
To this end, we observe that \eqref{convmut} implies the weak star convergence in $L^\infty(0,T)$ of the mean values, so~that
\Beq
  \norma\mubar_\infty \leq c \,.
  \label{stimamubar}
\Eeq
On the other hand, \eqref{intprimaz} and the regularity properties already established yield
\begin{align}
  & \< \mu , -\Delta v >_{\L2W}
  = - \< \dt\phi , v >_{\L2V}
  + \intQ (g-\sigma\phi) v
  \non
  \\
  & \leq c \, \norma v_{\L2V}
  \quad \hbox{for every $v\in\L2Z$}.
  \label{perlemma} 
\end{align}
Now, for a given $v\in\L2W$, we solve the (time dependent) elliptic problem
\Beq
  - \Delta z = v - \vbar
  \quad \aeQ \,, \quad 
  \zbar = 0
  \quad \aet
  \aand
  \dn z = 0 
  \quad \hbox{on $\Sigma$}.
  \non
\Eeq
This problem has a unique solution $z \pier{{}= \calN ( v - \vbar) }\in\L2W$ \pier{(cf.~\eqref{neumann} and \eqref{defN}).}
Moreover, the equation itself yields $\Delta z\in\L2W$, so that $z\in\L2Z$.
On the other hand
\Beq
  \ioT \< \mu(t) , v(t) > \, dt
  = \ioT \< \mu(t) , \vbar(t) > \, dt
  + \ioT \< \mu(t) , - \Delta z(t) > \, dt
  \non
\Eeq
and we easily estimate the terms on the \rhs.
From \eqref{stimamubar} we derive~that
\Beq
  \ioT \< \mu(t) , \vbar(t) > \, dt
  = \ioT |\Omega| \, \mubar(t) \, \vbar(t) \, dt
  \leq c \, \norma\vbar_{L^2(0,T)}
  \leq c \, \norma v_{\L2\Vp}
  \non
\Eeq
while \eqref{perlemma} immediately yields
\Beq
  \ioT \< \mu(t) , - \Delta z(t) > \, dt
  \leq c \, \norma z_{\L2V}
  \leq c \, \norma{v-\vbar}_{\L2\Vp}
  \leq c \, \norma v_{\L2\Vp} \,.
  \non
\Eeq
Hence, we have that
\Beq
  \ioT \< \mu(t) , v(t) > \, dt
  \leq c \, \norma v_{\L2\Vp} \,.
  \non
\Eeq
Since $v\in\L2W$ is arbitrary, we can apply Lemma~\ref{Permu} 
and conclude that $\mu$ belongs to~$\L2V$.
It remains to recover the regularity \eqref{regwz} of $w$ and~\eqref{secondaz}.
To this end, we rewrite \eqref{secondazbis} in the~form
\begin{align}
  &\pier{ \iO w \, (-\Delta v + v) = \iO \bigl( \mu - \beta'(\phi) \, w + (\lambda-\nu+1) w \bigr) \, v }
  \non
  \\
  & 
  \pier{\quad  \hbox{ for every $v\in W$ and \aet} \,. }
  \non
\end{align}
Then, we can apply Lemma~\ref{Perw} and infer that
\begin{align}
  & w(t) \in W \,, \quad 
  - \Delta w(t) + w(t)
  = (\mu-\beta'(\phi)\,w+(\lambda-\nu+1)w)(t)
  \non
  \\
  & \aand \normaW{w(t)} \leq \CO \, \norma{(\mu-\beta'(\phi)\,w+(\lambda-\nu+1)w)(t)}
  \quad \aat \,.
  \non
\end{align}
Then, $w \in\L2W$ and \eqref{secondaz} is satisfied.
This completes the proof of Theorem~\ref{Asymptotics}.

\pier{\subsection{Proof of Theorem~\ref{ErrEst}}}%

This subsection is devoted to the proof of the error estimate~\eqref{error}.
Then, for $\tau \in (0,1) $ let $\soluzt$ be a solution to Problem \Pbl\
that \gianni{satisfies} the estimates \Stability, and let  \pier{$\soluz$} denote the unique solution to Problem \Pblz. 
\gianni{In the sequel, we often account for the \Lip\ continuity of the nonlinearities on the range of $\phit$ and~$\phi$
since these functions are uniformly bounded due to \eqref{stab} and~\eqref{regphiz}.
To this concern, we point out once and for all that the corresponding \Lip\ constants are estimated 
by a constant that depends only on $\Omega$, $T$, the structure of the system and the norms of the data corresponding to~\eqref{hpdati}.
A~similar remark holds for some uniform bounds, like the $L^\infty$ norm of~$\beta(\phit)$.
To perform our proof of \eqref{error}, we set for convenience}
\Beq
  \overline{\phi} := \phi_\tau - \phi\,, \quad 
  \overline{\mu} := \mu_\tau - \mu \, , \quad
  \overline{w} := w_\tau  - w 
  \non
\Eeq
and take the difference between \eqref{prima}\gianni{, written for $\soluzt$ with $v\in Z$,} and 
\eqref{primaz}, {finding}~that
\Beq
  \tau \< \dt^2 \phit , v > + \< \dt\overline{\phi} , v >
  + \sigma \iO \overline{\phi}\,  v
= \< \overline{\mu} , \Delta v> 
  \quad \hbox{\hbox{for every $v\in Z$ and \aet}}\,. 
  \label{diffprima}
\Eeq
First, we take $v= 1/|\Omega|$, which is feasible since the constant functions
are in~$Z$, and we obtain the identity
\Beq
\tau (\phit)''_\Omega + \overline{\phi}^{\, \prime}_\Omega  + \sigma \, \overline{\phi}_\Omega
  = 0
    \label{diffprimabar}
\Eeq
for the mean value~$\overline{\phi}_\Omega$. 
Testing by $\overline{\phi}_\Omega$
and integrating over $(0,t)$, $t\in (0,T]$, 
from the initial conditions \gianni{$(\phit)'_\Omega(0)=(\rhoz)_\Omega $ and $\overline{\phi}_\Omega (0)=0$ (cf.\ \eqref{cauchy} and \eqref{cauchyz})}
and Young's inequality, we arrive~at 
\begin{align}
 & \frac12 \, | \overline{\phi}_\Omega (t)|^2 
  + \sigma\iot |\overline{\phi}_\Omega|^2 
  =  - \tau (\phit)'_\Omega  (t) \, \overline{\phi}_\Omega (t)  
  +  \tau \iot   (\phit)'_\Omega  \overline{\phi}^{\, \prime}_\Omega
 \non
 \\ 
 & \leq  \frac14 \, | \overline{\phi}_\Omega (t)|^2 + 
   c \, \tau^2 \, \norma{\dt \phit }_{\L\infty{\Vp}}^2 \, 
  + c\, \tau \, \norma{\dt \phit}_{\L2{\Vp}} \, \norma{\dt \overline{\phi}}_{\L2{\Vp}}
  \non
\end{align}
\aat, where we have used the norms in $\Vp$ for the mean values, which are constant in space. 
Hence, from \eqref{stab}--\eqref{regphiz} it follows~that 
\Beq
  \norma{\overline{\phi}_\Omega}_{C^0([0,T])}^2
  + \norma{\overline{\phi}_\Omega}^2_{L^2 (0,T)}\, 
    \leq \, c \, \tau
    \,.
  \label{1diffstima}
\Eeq
Next, we multiply \eqref{diffprimabar} by $v\in Z$ and integrate over~$\Omega$\gianni{; then we} 
subtract it to \eqref{diffprima}. 
Afterwards, we choose  
$v = \calN(\overline{\phi} - \overline{\phi}_\Omega)$ 
\gianni{in this difference}
and point out that this test function is in $\L\infty Z$ since 
$\overline{\phi} - \overline{\phi}_\Omega \in \L\infty W$ by 
\eqref{regphi} and \eqref{regphiz}. 
Besides, $ \calN(\overline{\phi} - \overline{\phi}_\Omega)$ has zero mean value due to~\eqref{defN}. 
Hence, recalling the properties \PropN, we have \aet~that
\Beq
  \frac 12 \, \ddt \, \normaVp{\overline{\phi} - \overline{\phi}_\Omega}^2
  + \sigma \, \normaVp{\overline{\phi} - \overline{\phi}_\Omega}^2
  = - \< \overline{\mu} , \overline{\phi} - 
  \overline{\phi}_\Omega> 
  - \tau \<\dt^2 \phit  , \calN 
  (\overline{\phi} - \overline{\phi}_\Omega) >  \,.
  \label{testdiffprima}
\Eeq
At the same time, we write \eqref{quarta} for both $\phit , \mut $ and $\phi, \mu $; then, we test the difference by $ \overline{\phi} - \overline{\phi}_\Omega $.
{Upon} rearranging, we obtain~that
\begin{align}
  & \iO |\Delta\overline{\phi}|^2
   = - \iO \nabla \bigl(
      \beta(\phit) - \beta(\phi)
    \bigr) \cdot \nabla\overline{\phi}
  \non
  \\
  & \quad {}
  + \iO \bigl(
    \beta'(\phit) \Delta\phit
    - \beta'(\phi) \Delta\phi
  \bigr) (\overline{\phi} - \overline{\phi}_\Omega)
  + (2\lambda-\nu) \iO |\nabla\overline{\phi}|^2 
  \non
  \\  
  & \quad {}
  {{}- \iO \bigl(
    \psi(\phit)
    - \psi(\phi)
  \bigr) (\overline{\phi} - \overline{\phi}_\Omega)}
  + \< \overline{\mu} , \overline{\phi} - 
  \overline{\phi}_\Omega>
  \label{testdiffsecondabis}
\end{align}
where the function $\psi$ is defined by 
\Beq
  \psi(s) := \beta(s) \beta'(s) - \lambda s \beta'(s) + (\nu-\lambda) \beta(s) + (\lambda^2-\lambda\nu) s
  \quad \hbox{for {every}  $s\in\erre$.}
  \label{defpsi}
\Eeq
Then, adding \eqref{testdiffprima} and \eqref{testdiffsecondabis} to each other 
leads to the cancellation of the two terms involving~$\overline{\mu}$. 
Moreover, the integration with respect
to time with the help of initial conditions~\eqref{cauchy} and \eqref{cauchyz}~gives 
\begin{align}
  &\frac 12 \, \normaVp{(\overline{\phi} - \overline{\phi}_\Omega)(t)  }^2
  + \sigma \iot \normaVp{\overline{\phi} - \overline{\phi}_\Omega}^2
  +  \intQt |\Delta\overline{\phi}|^2  
  \non\\
  &{}= 
  - \tau \<\dt \phit (t)  , \calN 
  (\overline{\phi} - \overline{\phi}_\Omega)(t) >  
+ \tau \iot  \<\dt \phit  , \dt \calN 
  (\overline{\phi} - \overline{\phi}_\Omega) > 
  \non\\
  &\quad {}   - \intQt \nabla \bigl(
      \beta(\phit) - \beta(\phi)
    \bigr) \cdot \nabla\overline{\phi}
 \, + \intQt \bigl(
    \beta'(\phit) \Delta\phit
    - \beta'(\phi) \Delta\phi
  \bigr) (\overline{\phi} - \overline{\phi}_\Omega)
  \non
  \\
  & \quad {} 
  + (2\lambda-\nu) \intQt |\nabla\overline{\phi}|^2 
  {}- \intQt \bigl(
    \psi(\phit)
    - \psi(\phi)
  \bigr) (\overline{\phi} - \overline{\phi}_\Omega)\, .
  \label{testdiffterza}
\end{align}  
Now, we can estimate the integrals on the \rhs\ of \eqref{testdiffterza} taking into account that 
$\phit$ and~$\phi$ are uniformly bounded (see \eqref{stab} and \eqref{convphit}), 
so that we can make use of the \Lip\ continuity of the nonlinearities on their range.
To begin with, by Young's inequality, \eqref{normaVp} and \eqref{stabbis} we have that
\begin{align}
 & - \tau \<\dt \phit (t)  , \calN (\overline{\phi} - \overline{\phi}_\Omega)(t) >  
 \le c \, \tau \, \|\dt \phit (t)\|_{\Vp} \|(\overline{\phi} - \overline{\phi}_\Omega)(t)\|_{*}
 \non \\
&{} \le \frac14 \, \|(\overline{\phi} - \overline{\phi}_\Omega)(t)\|_{*}^2
 +  c\, \tau^2 \, \|\dt \phit \|^2_{\L\infty {\Vp}}
 \le \frac14 \, \|(\overline{\phi} - \overline{\phi}_\Omega)(t)\|_{*}^2
 +  c\, \tau \, .
 \label{pier1}
\end{align}
On the other hand, using \eqref{stab} and \eqref{regphiz} we deduce that
\begin{align}
& \tau \iot  \<\dt \phit  , \dt \calN (\overline{\phi} - \overline{\phi}_\Omega) > 
\,=\, \tau \iot  \<\dt \phit  , \calN 
  (\dt \overline{\phi} - \overline{\phi}^{\, \prime}_\Omega) >  \non \\
&\le\, c\, \tau \iot  \| \dt \phit \|_{\Vp} 
  \norma{\dt \overline{\phi} - \overline{\phi}^{\, \prime}_\Omega}_* \,\le\, c\, \tau  \| \dt \phit \|_{\L2{\Vp}} 
  \norma{\dt \overline{\phi}}_{\L2{\Vp}}\,\le\, c\, \tau .\label{pier2}
\end{align}
To deal with the next term, we observe that $\beta'$ is nonnegative and recall
the continuity of the embedding $W\emb\Wx{1,4}$ and the bound \eqref{stab} for~$\phit$.
In fact, by virtue of the \Holder\ inequality and the compactness inequality~\eqref{compact}, for every $\delta>0$ it holds~that
\begin{align}
  & - \intQt \nabla \bigl(
    \beta(\phit) - \beta(\phi)
  \bigr) \cdot \nabla\overline{\phi}
  = - \intQt \bigl(
    \beta'(\phit)\nabla\phit
    - \beta'(\phi)\nabla\phi
  \bigr) \cdot \nabla\overline{\phi}
  \non
  \\
  & = - \intQt \bigl(
    \beta'(\phit) - \beta'(\phi)
  \bigr) \nabla\phit \cdot \nabla\overline{\phi}
  - \intQt \beta'(\phi) |\nabla\overline{\phi}|^2
  \leq c \intQt |\gianni{\overline\phi}| \, |\nabla\phit| \, |\nabla\overline{\phi}|
  \non
  \\
  & \leq c \, \iot \norma{\overline{\phi}}_4 \, \norma{\nabla\phit}_4 \, \norma{\nabla\overline{\phi}}
  \gianni{
  {}\leq c \iot \normaW\phit \, \normaV{\overline\phi}^2
  }
  \non
  \\
  & \leq c \iot \normaV{\overline{\phi}}^2
  \leq \delta \intQt |\Delta\overline{\phi}|^2
  + \cdelta \iot \normaVp{\overline{\phi}}^2 
  \non
\end{align}
\gianni{and the last integral can be estimated using \eqref{1diffstima} as follows}
\begin{align}
&\iot \normaVp{\overline{\phi}}^2 
  \leq 2 \iot \normaVp{\overline{\phi} - \overline{\phi}_\Omega}^2 
  + 2 \iot \normaVp{\overline{\phi}_\Omega}^2 \non \\
  &\leq 2 \iot \normaVp{\overline{\phi} - \overline{\phi}_\Omega}^2 + c \ioT |\overline{\phi}_\Omega|^2  
\leq 2 \iot \normaVp{\overline{\phi} - \overline{\phi}_\Omega}^2 + c\, \tau \, .
\non
\end {align}
Therefore, it turns out that
\Beq
- \intQt \nabla \bigl(
    \beta(\phit) - \beta(\phi)
  \bigr) \cdot \nabla\overline{\phi}
  \leq \delta \intQt |\Delta\overline{\phi}|^2
  + \cdelta \iot \normaVp{\overline{\phi} - \overline{\phi}_\Omega}^2
  + \cdelta \, \tau \, .
  \label{pier3}
\Eeq
By proceeding and arguing as above we discuss the term
\begin{align}
  & \intQt \bigl(
    \beta'(\phit) \Delta\phit
    - \beta'(\phi) \Delta\phi
  \bigr) (\overline{\phi} - \overline{\phi}_\Omega)
  \non
  \\
  & = \intQt \bigl( \beta'(\phit) - \beta'(\phi) \bigr) \Delta\phit \, (\overline{\phi} - \overline{\phi}_\Omega)
  + \intQt \beta'(\phi) \Delta\overline{\phi} \, (\overline{\phi} - \overline{\phi}_\Omega)
  \non
  \\
  & \leq c \iot \norma{\overline{\phi}}_4 \, \norma{\Delta\phit} \, \norma{\overline{\phi} - \overline{\phi}_\Omega}_4
  + c \iot \norma{\Delta{\overline{\phi}}} \, \norma{\overline{\phi} - \overline{\phi}_\Omega}
  \non
  \\
  & \leq c \iot \normaV{\overline{\phi}}\, \normaV{\overline{\phi} - \overline{\phi}_\Omega}
  + c \iot \norma{\Delta{\overline{\phi}}} \, \norma{\overline{\phi} - \overline{\phi}_\Omega}
  \non
\end{align} 
whence, by Young's inequality\gianni{, \eqref{compact} and~\eqref{1diffstima}}, 
\begin{align}
  & \intQt \bigl(
    \beta'(\phit) \Delta\phit
    - \beta'(\phi) \Delta\phi
  \bigr) (\overline{\phi} - \overline{\phi}_\Omega)
  \non
  \\
  \separa
  & \leq c \iot \normaV{\overline{\phi}}^2  + 
  c \iot \normaV{\overline{\phi} - \overline{\phi}_\Omega}^2
  + \delta \intQt |\Delta\overline{\phi}|^2
  + \cdelta \iot \norma{\overline{\phi} - \overline{\phi}_\Omega}^2
  \non
  \\
  \separa
  & \leq 4 \delta \intQt |\Delta\overline{\phi}|^2 
  + \cdelta  \iot \normaVp{\overline{\phi}}^2  
  + \cdelta \iot \normaVp{\overline{\phi} - \overline{\phi}_\Omega}^2. 
  \non
    \\
  & \leq 4 \delta \intQt |\Delta\overline{\phi}|^2
  + \cdelta \iot \normaVp{\overline{\phi} - \overline{\phi}_\Omega}^2
  +\cdelta \, \tau \,. 
  \label{pier4}
  \end{align}
As $\psi$ is \Lip\ continuous on the range of the solutions {under consideration}, 
we can \gianni{similarly} handle the last two terms of \eqref{testdiffterza} as follows
\begin{align}
  &(2\lambda-\nu) \intQt |\nabla\overline{\phi}|^2 
  {}- \intQt \bigl(\psi(\phit) - \psi(\phi)
  \bigr) (\overline{\phi} - \overline{\phi}_\Omega)
  \non
  \\
  & \leq c \iot \normaV{\overline{\phi}}^2 \, 
  + c \iot \norma{\overline{\phi}} \, \norma{\overline{\phi} - \overline{\phi}_\Omega}
 \leq c \iot \normaV{\overline{\phi}}^2 \, 
  + c \iot \norma{\overline{\phi} - \overline{\phi}_\Omega}^2
  \non
    \\
  & \leq 2 \delta \intQt |\Delta\overline{\phi}|^2 + \cdelta  \iot \normaVp{\overline{\phi}}^2 
   + \cdelta \iot \normaVp{\overline{\phi} - \overline{\phi}_\Omega}^2
  \non
    \\
  & \leq 2 \delta \intQt |\Delta\overline{\phi}|^2
  + \cdelta \iot \normaVp{\overline{\phi} - \overline{\phi}_\Omega}^2 + \cdelta \, \tau \,. 
  \label{pier5}
  \end{align}
Then, we collect all of the estimates~\eqref{pier1}--\eqref{pier5} to control the \rhs\ of \eqref{testdiffterza}\gianni{; then we} 
choose $\delta$ small enough (e.g., $\delta= 1/14$) and apply the Gronwall lemma to conclude~that
\Beq
  \norma{\overline{\phi} - \overline{\phi}_\Omega}_{\C0\Vp}
  + \norma{\Delta\overline{\phi}}_{\L2H}
  \leq c \, \tau^{1/2} .
  \label{provvisoria}
\Eeq
By combining this estimate with \eqref{1diffstima} and applying the elliptic regularity theory, we find out that
\Beq
  \norma{\overline{\phi}}_{\C0\Vp {\cap}\L2W}
  \leq c \, \tau^{1/2} . 
  \label{pier6}
\Eeq
Next, we recall \eqref{terza} for $\phit, \wt$ and \eqref{terzaz} for $\phi, w $
and take their difference obtaining
\Beq
\overline{w} =  - \Delta\overline{\phi} + \beta(\phit) - \beta (\phi)- \lambda\overline{\phi}
  \quad \aeQ . \non 
\Eeq
Then, \gianni{the Lipschitz continuity of~$\beta$ and \eqref{pier6} yield} the additional estimate 
\Beq
  \norma{\overline{w}}_{\L2H}
  \leq c \, \tau^{1/2} . 
  \label{pier7}
\Eeq
Finally, arguing with \gianni{\eqref{seconda} for $\soluzt$ and \eqref{secondaz} for $\soluz$,}
it is easy to derive the equality
\begin{align}
 &\ioT \< \overline \mu , v > =  \intQ \overline w \, (- \Delta v + (\nu-\lambda)v)   + \intQ (\beta'(\phit) - \beta'(\phi)) \wt v
 + \intQ \beta'(\phi) \overline w \,  v 
  \label{pier8}
\end{align} 
for every $v\in \L2W$. 
Concerning the \rhs\ of \eqref{pier8}, with the help of H\"older's inequality, the continuous embedding $W\hookrightarrow L^\infty (\Omega) $, the bound in \eqref{stab} and the estimates \eqref{pier6} and \eqref{pier7}  we infer that
\begin{align}
 &\biggl| \intQ \overline w \, (- \Delta v + (\nu-\lambda)v)   + \intQ (\beta'(\phit) - \beta'(\phi)) \wt v
 + \intQ \beta'(\phi) \overline w  v \biggr| \non \\
 &\leq c \, \norma{\overline w}_{\L2H} \norma{v}_{\L2W} +
 c \, \norma{\wt}_{\L\infty H} \ioT \norma{\overline \phi}_\infty \norma{v} 
 \non\\
 &\quad{}+ 
    \norma{\beta' (\phi)}_{C^0(\overline{Q})} \iot \norma{\overline w} \norma{v} 
\non \\
 &\leq c\, \tau^{1/2} \norma{v}_{\L2W} +
 c \Bigl( \norma{\overline \phi}_{\L2W} +
 \norma{\overline w}_{\L2H} \Bigr) \norma{v}_{\L2H} \non \\
 &\leq c \,\tau^{1/2} \norma{v}_{\L2W} +
 c \, \tau^{1/2} \norma{v}_{\L2H} . 
  \label{pier9}
\end{align} 
Therefore, from \eqref{pier8} and \eqref{pier9} it is straightforward to show that 
\Beq
  \norma{\overline{\mu}}_{\L2{W^*}}
  \leq c \, \tau^{1/2}  
  \label{pier10}
\Eeq
which, along with \eqref{pier6} and \eqref{pier7}, ends the proof of 
\eqref{error} and Theorem~\ref{ErrEst}.%


\section*{Acknowledgments}

\pier{PC acknowledges the support of the Next Generation EU Project No.P2022Z7ZAJ 
(A unitary mathematical framework for modelling muscular 
dystrophies), the RISM (Research Institute for Mathematical Sciences,
an International Joint Usage/Research Center located in Kyoto  University)
and the GNAMPA (Gruppo Nazionale per l'Analisi Matematica, la Probabilit\`{a} e le loro 
Applicazioni) of INdAM (Istituto Nazionale di Alta Matematica).}


\End{document}

 
\Beq
  \hbox{$\beta$ and $\beta''$ are monotone}.
  \label{betamonotona}
\Eeq
Moreover, $\Beta$ grows faster than $|s|^3$ as $|s|$ tends to infinity,
so that 
\Beq
  \hbox{$F$ is bounded from below}.
  \label{bddbelow}
\Eeq